\documentclass[12pt]{amsart}

\usepackage{amssymb, verbatim}

\title[On factorization of $q$-d.e. for continuous $q$-ultraspherical polynomials]
{On factorization of $q$-difference equation for continuous $q$-ultraspherical polynomials}

\author[Area]{I. Area}
\author[Atakishiyeva]{M.K. Atakishiyeva}
\author[Rodal]{J. Rodal}

\address[Area and Rodal]{Departamento de Matem\'atica Aplicada II, E.T.S.E. Telecomunicaci\'on, Universidade de Vigo, 36310--Vigo, Spain.}
\email[Area]{area@dma.uvigo.es}
\email[Rodal]{jrodal@edu.xunta.es}

\address[Atakishiyeva]{Facultad de Ciencias, Universidad Aut\'onoma del Estado de Morelos, C.P. 62250 Cuernavaca, Morelos, M\'exico.}
\email[Atakishiyeva]{mesuma@servm.fc.uaem.mx}



\begin{document}

\begin{abstract}

We prove that a customary Sturm-Liouville form of second-order $q$-dif\-ference
equation for the continuous $q$-ultraspherical polynomials $C_n(x;\beta|\,q)$
of Rogers can be written in a factorized form in terms of some explicitly
defined $q$-difference operator ${\mathcal D}_x^{\beta,\,q}$. This reveals the
fact that the continuous $q$-ultraspherical polynomials $C_n(x;\beta|\,q)$
are actually governed by the $q$-difference equation ${\mathcal D}_x^{\beta,\,q}\,
C_n(x;\beta|\,q)= \left(q^{-n/2}+\beta\,q^{n/2}\right)\,C_n(x;\beta|\,q)$,
which can be regarded as a square root of the equation, obtained from its
original form.

\end{abstract}


\maketitle

\setcounter{equation}{0}

\section{Introduction}

\medskip

It is well known that for many purposes it proves practical, as in
the case of linear second-order ordinary differential equations, to
represent the difference equation of hypergeometric type for classical
orthogonal polynomials in Sturm-Liouville (or self-adjoint) form \cite{NU}
$$
\frac{\Delta}{\Delta x(s-1/2)}\left[ \sigma(s)\,\rho(s)\,\frac{\nabla\,f(s)}
{\nabla\,x(s)}\,\right]\,+\,\lambda\,\rho(s)\,f(s)\,=0\,,    \eqno(1.1)
$$
where $\Delta\,y(s):=y(s+1)-y(s)$ and $\nabla\,y(s):=y(s)-y(s-1)$
(we employ standard notations of the theory of special functions, see,
for example, \cite{GR} or \cite{AAR}).

The important feature of this form (1.1) is that it requires the introduction
of a function $\rho(s)$ through the Pearson-type difference equation
$$
\frac{\Delta}{\Delta x(s-1/2)}\left[\,\sigma(s)\,\rho(s)\,\right]\,=
\,\tau(s)\,\rho(s)\,,                                           \eqno(1.2)
$$
with polynomials $\sigma(s)$ and $\tau(s)$ of respective degrees at most two
and one, which cha\-racterize an original form of the difference equation (1.1).
The full importance of the self-adjoint form (1.1) becomes apparent when one
takes into account that the same function $\rho(s)$ enables one to formulate
the orthogonality property of solutions of equation (1.1). Moreover, one can
construct explicit representation (\cite{NU}, p.66)
$$
f_n(s):= \frac{B_n}{\rho(s)}\,\frac{\nabla}{\nabla\,x_1(s)}\,\cdots\,
\frac{\nabla}{\nabla\,x_{n-1}(s)}\,\frac{\nabla}{\nabla\,x_n(s)}\,
\left[\,\rho(s+n)\,\prod_{k=1}^{n}\,\sigma(s+k)\,\right]
$$
in terms of the function $\rho(s)$ for the polynomial solutions $f_n(s)$ of
equation (1.1), which correspond to the values $\lambda_n:=-n\,\tau'-n(n-1)\,
\sigma^{''}/2$ of the parameter $\lambda$ (for a more detailed discussion of
this topic, see \cite{NU}).

An example to illustrate this point is provided by the continuous $q$-Hermite
polynomials of Rogers,
$$
H_n(x|\,q):= \sum_{k=0}^n\,\left[\!\!\begin{array} {c} n \\ k \end{array}
\!\!\right]_q e^{{\rm i}(n-2k)\,\theta}\,, \quad\quad 0<q<1\,,      \eqno(1.3)
$$
which are orthogonal on the finite interval $-1\leq x := \cos \theta \leq 1$
with respect to the weight function
$$
\widetilde w (x|\,q):= \,\frac{1}{\sin \theta}\,
\left(\,e^{2{\rm i}\theta},\,e^{-2{\rm i}\theta}; q \right)_\infty \, . \eqno(1.4)
$$
These polynomials $H_n(x|\,q)$ satisfy the following $q$-difference equation
$$
D_{q}\,\left[{\widetilde w}(x|\,q)\,D_q\,H_n(x|\,q)\right] =
\frac{4\,q\,(1-q^{-n})}{(1-q)^2}\,H_n(x|\,q)\,{\widetilde w}(x|\,q)\,, \eqno(1.5)
$$
written in self-adjoint form (1.1) (see \cite{KS}, p.115). The symbol $D_q$ 
in (1.5) is the conventional notation for the Askey-Wilson divided-difference 
operator (see, for example, \cite{AAR}, p.529), defined as
$$
D_{q}\,  f(x) := \frac{\delta_{q}\,f(x)}{\delta_{q}\,x}\,,\quad\quad
\delta_q\,g(e^{{\rm i}\,\theta}):= g(q^{1/2}\,e^{{\rm i}\,\theta}) -
g(q^{-1/2}\,e^{{\rm i}\,\theta})\,,\quad x=\cos\theta \,.            \eqno(1.6)
$$

As was observed in \cite{AK}, one may eliminate the weight function
$\widetilde w(x|\,q)$ from (1.5) by utilizing its readily verified
property that
$$
\exp \left(\pm \,{\rm i}\ln q^{1/2}\,\partial_{\theta}\right){\widetilde w}(x|\,q)
= - \frac{e^{\pm 2{\rm i}\theta}}{\sqrt{q}}\,{\widetilde w}(x|\,q)\,.    \eqno (1.7)
$$
It should be noted that following \cite{AK} we find it more convenient to write
(1.7) (and subsequent $q$-difference equations) in terms of the shift operators
(or the operators of the finite displacement, \cite{LL}) $e^{\pm a\,\partial_{\theta}}
\,g(\theta):= g(\theta \pm a)$ with respect to the variable $\theta$.

This elimination of the weight function ${\widetilde w}(x|\,q)$ from (1.5) yields
the following $q$-difference equation
$$
\frac{1}{2{\rm i}\sin\theta}\,\left[\,\frac{e^{{\rm i}\theta}}
{1-q\,e^{-2{\rm i}\theta}}\,\left(e^{{\rm i}\ln q\,\partial_{\theta}}\,
- 1\right)+ \,\frac{e^{-{\rm i}\theta}}{1 - q\,e^{2{\rm i}\theta}}\,
\left(1 - e^{-{\rm i}\ln q\,\partial_{\theta}}\,\right)\right]\,H_n(x|\,q)$$
$$ = \left( q^{-n}-1\right)\,H_n(x|\,q)                                 \eqno(1.8)
$$
for the continuous $q$-Hermite polynomials $H_n(x|\,q)$. The resultant
$q$-difference equation (1.8) then admits factorization of the form
$$
\left({\mathcal D}_x^{\,q}\right)^2 H_n(x|\,q) = q^{-n}\,H_n(x|\,q)\,,      \eqno(1.9)
$$
where the $q$-difference operator ${\mathcal D}_x^{\,q}$ is equal to
$$
{\mathcal D}_x^{\,q} := \,\frac{1}{1- e^{-2{\rm i}\theta}}\,\,e^{\,{\rm i}\,
\ln q^{1/2}\,\partial_{\theta}} + \frac{1}{1- e^{2{\rm i}\theta}}\,\,
e^{-{\rm i}\,\ln q^{1/2}\,\partial_{\theta}}                      $$
$$
\equiv \frac{1}{2{\rm i}\sin\theta}\,\left(e^{{\rm i}\theta}\,
e^{{\rm i}\,\ln q^{1/2}\,\partial_{\theta}} - e^{{-\rm i}\theta}\,
e^{{-\rm i}\,\ln q^{1/2}\,\partial_{\theta}}\,\right)\,, \quad
x=\cos\theta\,.                                                        \eqno(1.10)
$$
This means that the continuous $q$-Hermite polynomials are in fact
governed by a simpler $q$-difference equation,
$$
{\mathcal D}_x^{\,q} H_n(x|\,q) = q^{-n/2}\,H_n(x|\,q)\,,                  \eqno(1.11)
$$
which represents a "square root" of (1.8) or (1.9).

This curious interrelation between two $q$-difference equations (1.5) and (1.11),
stu\-died in detail in \cite{AK}, leads to the natural question whether the continuous
$q$-Hermite polynomials $H_n(x|\,q)$ represent the exceptional case or there
exist other instances of orthogonal polynomials from the Askey $q$-scheme
\cite{KS}, which admit the same type of factorization in corresponding
$q$-difference equations for them.

The present paper is aimed at proving that the continuous $q$-ultraspherical
(Rogers) polynomials $C_n(x;\beta|\,q)$ exhibit the same property of factorization
as the continuous $q$-Hermite polynomials $H_n(x|\,q)$. The next section collects
those known facts about the $q$-ultraspherical  polynomials $C_n(x;\beta|\,q)$ and
their $q\to1$ limit counterpart, the Gegenbauer (ultraspherical) polynomials
$C_n^{(\gamma)}(x)$, which are needed in section 3 for proving that a $q$-difference
equation for the $C_n(x;\beta|\,q)$, derived from its appropriate self-adjoint form
like (1.1), does admit a factorization of the type (1.9). In the concluding section
4 we briefly discuss some special and limit cases of the parameter $\beta$, which
are related with other well-known families of $q$-polynomials.
\medskip

\setcounter{equation}{0}
\section{Rogers and Gegenbauer polynomials}
\medskip

To proceed further we need to recall in this section some standard facts
about continuous $q$-ultraspherical (Rogers) polynomials and their $q\to1$
limit counterpart, Gegenbauer (ultraspherical) polynomials. The continuous
$q$-ultraspherical polynomials
$$
C_n(x;\beta|\,q):= \sum_{k=0}^n\,\frac{(\beta;q)_k\,(\beta;q)_{n-k}}
{(q;q)_k\,(q;q)_{n-k}}\,e^{{\rm i}(n-2k)\,\theta}\,, \quad\quad 0<q<1\,,\eqno(2.1)
$$
are known to be orthogonal on the finite interval $-1\leq x := \cos \theta
\leq 1$,
$$
\frac{1}{2\pi}\,\int_{-1}^{1}\,C_m(x;\beta|\,q)\,C_n(x;\beta|\,q)\,
{\widetilde w} (x;\beta|\,q)\,dx\,=\,d_n^{-1}(\beta;q)\, \delta_{mn}\,,$$
$$ d_n(\beta;q):= \frac{(1-\beta q^n)}{(1-\beta)}\,\frac{(q;q)_n}{(\beta^2;q)_n}\,
\frac{(\beta^2,q;q)_{\infty}}{(\beta,\beta q;q)_{\infty}}\,,
 \quad\quad |\beta| < 1 \,,                         \eqno(2.2)
$$
with respect to the weight function (see, for example, \cite{KS}, p.86)
$$
\widetilde w (x;\beta|\,q):= \,\frac{1}{\sin \theta}\,
\frac {\left(\,e^{2{\rm i}\theta},\,e^{-2{\rm i}\theta};q \right)_\infty}
{\left(\,\beta\,e^{2{\rm i}\theta},\,\beta \,
e^{-2{\rm i}\theta}; q \right)_\infty}\,\,.                       \eqno(2.3)
$$
They satisfy the Sturm-Liouville type $q$-difference equation
$$
D_{q}\,\left[{\widetilde w}(x;\beta q|\,q)\,D_q\,C_n(x;\beta|\,q)\right] =
\lambda_n(\beta)\,C_n(x;\beta|\,q)\,{\widetilde w}(x;\beta|\,q)   \eqno(2.4)
$$
with eigenvalues $\lambda_n(\beta):= 4\,q\,(1-q^{-n})\,(1-\beta^2 q^n)/(1-q)^2$\,
(see, for example, \cite{KS}, p.86). Note that the $D_q$ in (2.4) is the same
Askey--Wilson divided-difference operator, defined above in (1.6), namely,
$$
D_{q}= \frac{\sqrt{q}}{{\rm i}(1-q)}\,\frac{1}{\sin\theta}\left( e^{\,{\rm i}
\ln q^{1/2}\,\partial_{\theta}}- e^{-{\rm i}\ln q^{1/2}\,\partial_{\theta}}\right),
\quad \quad \partial_{\theta}\equiv\frac{d}{d\theta}\,.                  \eqno(2.5)
$$
Observe also that one readily derives from definition (2.3) the relation
$$
{\widetilde w}(x;\beta q|\,q)= \left[(1+\beta)^2 - 4 \beta x^2\right]
{\widetilde w}(x;\beta|\,q)                                              \eqno(2.6)
$$
between the weight functions ${\widetilde w}(x;\beta|\,q)$  with the two distinct
parameters $\beta$ and $\beta q$. Therefore a $q$-analogue of the factor $\sigma(s)$
from the self-adjoint equation (1.1) in the case of the $q$-difference equation (2.4)
is just
$$
\sigma_q(x;\beta):= (1+\beta)^2 - 4\beta x^2 \,.
$$

If one sets $\beta=q^{\gamma}$ in (2.1) and then evaluates its limit as
$q\to 1$, this results in
$$
\lim_{q\to1 }\,C_n(x;q^{\gamma}|\,q)= C_n^{(\gamma)}(x)\,,
$$
where $C_n^{(\gamma)}(x)$ are the Gegenbauer polynomials:
$$
C_n^{(\gamma)}(x):= \sum_{k=0}^n\,\frac{(\gamma)_k\,(\gamma)_{n-k}}{k !\,(n-k)!}
\,\,e^{{\rm i}(n-2k)\,\theta}\,, \quad\quad x=\cos\theta\,.             \eqno(2.7)
$$
The self-adjoint form of the second-order differential equation for the
Gegenbauer polynomials (2.7) is known to be of the form
$$
\frac{d}{dx}\left[\,(1-x^2)\,w(x)\,\frac{d\,C_n^{(\gamma)}(x)}{dx}\right] +
n(n+2\gamma)\,w(x)\,C_n^{(\gamma)}(x)=0,                           \eqno(2.8)
$$
where $w(x):= (1-x^2)^{\gamma - 1/2}$, $\gamma > - 1/2$, is the orthogonality weight
function for the $C_n^{(\gamma)}(x)$ on the finite interval $-1\leq x \leq 1$. After
eliminating the weight function $w(x)$ from (2.8), one can rewrite it as
$$
\left[\,(1-x^2)\,\frac{d^2}{dx^2} - (2\,\gamma+1)\,x\,\frac{d}{dx}
+ n(n+2\gamma)\,\right]C_n^{(\gamma)}(x)=0\,.                      \eqno(2.9)
$$
In contrast to (2.8), this differential equation is evidently not self-adjoint;
but to transform it into self-adjoint equation (2.8) one needs only to multiply
it by $w(x)$ from the left and employ the readily verified identity
$$
w(x)\left[\,(1-x^2)\,\frac{d^2}{dx^2} - (2\,\gamma+1)\,x\,\frac{d}{dx}\,\right]\,
=\,\frac{d}{dx}\left[\,(1-x^2)\,w(x)\,\frac{d}{dx}\right]\,.             \eqno(2.10)
$$
\medskip

\setcounter{equation}{0}
\section{Factorization for Rogers polynomials}
\medskip

To eliminate the weight function $\widetilde w (x;\beta|\,q)$ from $q$-difference
equation (2.4), we employ first the relations
$$
\exp \left(\pm \,{\rm i}\ln q^{1/2}\,\partial_{\theta}\right)
{\widetilde w}(x;\beta q|\,q) = \frac {1}{\sqrt{q}}\,\left(1 - \beta\,
e^{\mp\, 2{\rm i}\theta}\right)\,\left(\beta q - e^{\pm \,2{\rm i}\theta}
\right)\,{\widetilde w}(x;\beta|\,q)\,,                       \eqno (3.1)
$$
which are straightforward to derive upon using the explicit expression (2.3) for
$\widetilde w (x;\beta|\,q)$ and relation (2.6). Substituting (3.1) into (2.4),
one obtains the $q$-difference equation
$$
\frac{1}{{\rm i}\sin\theta}\,\left[\,e^{{\rm i}\theta}\,\frac{\left(1- \beta\,
e^{-2{\rm i}\theta}\right)\,\left(1- \beta q\,e^{-2{\rm i}\theta}\right)}
{1-q\,e^{-2{\rm i}\theta}}\,\left(e^{{\rm i}\ln q\,\partial_{\theta}}\,-
1\right)\right.$$
$$
\left. +\, e^{-{\rm i}\theta}\,\frac{\left(1- \beta\, e^{2{\rm i}\theta}\right)
\,\left(1- \beta q\,e^{2{\rm i}\theta}\right)}{1-q\,e^{2{\rm i}\theta}}\,
\left(1 - e^{-{\rm i}\ln q\,\partial_{\theta}}\right)\right]\,C_n(x;\beta|\,q)$$
$$
= 2\, \left( q^{-n}-1\right)\,\left( 1- \beta^2\,q^{n}\right)\,C_n(x;\beta|\,q)\,\eqno(3.2)
$$
for the $q$-ultraspherical polynomials $C_n(x;\beta|\,q)$, which does not contain
the weight function $\widetilde w (x;\beta|\,q)$. This equation is a $q$-extension
of the second-order differential equation (2.9) for the Gegenbauer polynomials
$C_n^{(\gamma)}(x)$.

The next step is to use two simple trigonometric identities
$$
\frac{e^{\,\pm\,{\rm i}\theta}}{{\rm i}\sin\theta} =
\pm \,\frac{2}{1-e^{\,\mp\, 2{\rm i}\theta}}
$$
in order to write a $q$-difference operator on the left side of equation (3.2) as
$$
2\,\left[\frac{\left(1-\beta\,e^{-2{\rm i}\theta}\right)\left(1- \beta q\,
e^{-2{\rm i}\theta}\right)}{\left(1-\,e^{-2{\rm i}\theta}\right)\left(1-q\,
e^{-2{\rm i}\theta}\right)}\left(e^{{\rm i}\ln q\,\partial_{\theta}} - 1\right)
- \frac{\left(1- \beta e^{2{\rm i}\theta}\right)\left(1- \beta q e^{2{\rm i}
\theta}\right)}{\left(1- e^{2{\rm i}\theta}\right)\left( 1-q\,e^{2{\rm i}
\theta}\right)}\left(1 - e^{-{\rm i}\ln q\,\partial_{\theta}}\right)\right]$$
$$
= 2\,\left[\frac{\left(1-\beta\,e^{-2{\rm i}\theta}\right)\left(1- \beta\,
q\,e^{-2{\rm i}\theta}\right)}{\left(1-\,e^{-2{\rm i}\theta}\right)\left(1-q\,
e^{-2{\rm i}\theta}\right)}\,\,e^{{\rm i}\ln q\,\partial_{\theta}} + \frac{\left(1
- \beta\,e^{2{\rm i}\theta}\right)\left(1- \beta\,q\, e^{2{\rm i}\theta}\right)}
{\left(1- e^{2{\rm i}\theta}\right)\left( 1-q\,e^{2{\rm i}\theta}\right)}\,\,
e^{-{\rm i}\ln q\,\partial_{\theta}}\right.$$
$$
\left.- \frac{\left(1-\beta\,e^{-2{\rm i}\theta}\right)\left(1- \beta\, q\,
e^{-2{\rm i}\theta}\right)}{\left(1-\,e^{-2{\rm i}\theta}\right)\left(1-q\,
e^{-2{\rm i}\theta}\right)}- \frac{\left(1- \beta\, e^{2{\rm i}\theta}\right)
\left(1- \beta\,q\,e^{2{\rm i}\theta}\right)}{\left(1- e^{2{\rm i}\theta}\right)
\left( 1-q\,e^{2{\rm i}\theta}\right)}\right]\,.                   \eqno(3.3)$$
The last important step is to employ a readily verified identity
$$
\frac{1-\beta\,q\,e^{\mp\,2{\rm i}\theta}}{1- q\,e^{\mp\,2{\rm i}\theta}}
\,\,e^{\pm\,{\rm i}\,\ln q^{1/2}\,\partial_{\theta}}= e^{\pm\,{\rm i}\,
\ln q^{1/2}\,\partial_{\theta}}\,\,\frac{1-\beta\,e^{\mp\,2{\rm i}\theta}}
{1- e^{\mp\,2{\rm i}\theta}}                                    \eqno(3.4)
$$
for the shift operators $\exp\left({\pm\,{\rm i}\,\ln q^{1/2}\,\partial_{\theta}}
\right)$, which enter into first two terms in (3.3). With the aid of (3.4) one can
thus cast (3.3) into the form
$$
2\,\left[\frac{1-\beta\,e^{-2{\rm i}\theta}}{1-\,e^{-2{\rm i}\theta}}
\,e^{{\rm i}\ln q^{1/2}\,\partial_{\theta}}\, \frac{1-\beta\,e^{-2{\rm i}\theta}}
{1-\,e^{-2{\rm i}\theta}}\,e^{{\rm i}\ln q^{1/2}\,\partial_{\theta}}\,
+\,\frac{1-\beta\,e^{2{\rm i}\theta}}{1-\,e^{2{\rm i}\theta}}
\,e^{-{\rm i}\ln q^{1/2}\,\partial_{\theta}}\, \frac{1-\beta\,e^{2{\rm i}\theta}}
{1-\,e^{2{\rm i}\theta}}\,e^{-{\rm i}\ln q^{1/2}\,\partial_{\theta}}\right. $$
$$
\left. \, - \,\frac{\left(1-\beta\,e^{-2{\rm i}\theta}\right)\left(1- \beta q\,e^{-2{\rm i}
\theta}\right)}{\left(1-\,e^{-2{\rm i}\theta}\right)\left(1-q\,e^{-2{\rm i}\theta}\right)}
\,- \,\frac{\left(1-\beta\,e^{2{\rm i}\theta}\right)\left(1- \beta q\,e^{2{\rm i}\theta}\right)}
{\left(1-\,e^{2{\rm i}\theta}\right)\left(1-q\,e^{2{\rm i}\theta}\right)}\,\right] $$
$$
= 2\,\left[\frac{1-\beta\,e^{-2{\rm i}\theta}}{1-\,e^{-2{\rm i}\theta}}
\,e^{{\rm i}\ln q^{1/2}\,\partial_{\theta}}\, \frac{1-\beta\,e^{-2{\rm i}
\theta}}{1-\,e^{-2{\rm i}\theta}}\,e^{{\rm i}\ln q^{1/2}\,\partial_{\theta}}
\, - \frac{(1+q)(1-\beta)(\beta-q)}{(1+q)^2- 4 q x^2}\right. $$
$$
\left. + \,\frac{1-\beta\,e^{2{\rm i}\theta}}{1-\,e^{2{\rm i}\theta}}
\,e^{-{\rm i}\ln q^{1/2}\,\partial_{\theta}}\, \frac{1-\beta\,e^{2{\rm i}
\theta}}{1-\,e^{2{\rm i}\theta}}\,e^{-{\rm i}\ln q^{1/2}\,\partial_{\theta}}
- 1 - \beta^2 \right]\,.                                              \eqno(3.5)
$$
It is not hard to verify now that the above expression represents a product
of two $q$-difference operators,
$$
2\,\left[\left({\mathcal D}_x^{\beta,\,q}\,\right)^2 - (1+\beta)^2\right]=
2\,\left({\mathcal D}_x^{\beta,\,q} + 1 + \beta\right)\,
\left({\mathcal D}_x^{\beta,\,q} - 1 - \beta \right) \,,
$$
where ${\mathcal D}_x^{\beta,\,q}$ is equal to ({\it cf.} (1.10))
$$
{\mathcal D}_x^{\beta,\,q} := \frac{1-\beta\,e^{-2{\rm i}\theta}}{1-\,e^{-2{\rm i}\theta}}
\,e^{{\rm i}\ln q^{1/2}\,\partial_{\theta}} + \frac{1-\beta\,e^{2{\rm i}\theta}}
{1-\,e^{2{\rm i}\theta}}\,e^{-{\rm i}\ln q^{1/2}\,\partial_{\theta}}
\equiv {\mathcal D}_x^{q}\,+\,\beta\,{\mathcal D}_x^{1/q} $$
$$ \equiv \frac{1}{2{\rm i}\sin\theta}\,\left[ \left(e^{{\rm i}\theta} -
\beta\,e^{-{\rm i}\theta}\right)\,e^{{\rm i}\ln q^{1/2}\,\partial_{\theta}}-
\left(e^{-{\rm i}\theta} - \beta\,e^{{\rm i}\theta}\right)\,
e^{-{\rm i}\ln q^{1/2}\,\partial_{\theta}}\right]\,.     \eqno(3.6)
$$
Finally, taking into account that the factor $\left( q^{-n}-1\right)\left(1 - \beta^2
\,q^{n}\right)$ on the right side of (3.2) can be written as $\left( q^{-n/2}+ \beta\,
q^{n/2}\right)^2\,- \left( 1 + \beta \right)^2$, one arrives at the following factorized
form of equation (3.2):
$$
\left({\mathcal D}_x^{\beta,\,q}\right)^2 C_n(x;\beta|\,q) =  \left(q^{-n/2}
+ \beta\,q^{n/2}\right)^2 C_n(x;\beta|\,q)\,.                 \eqno(3.7)
$$

Note that the operator $({\mathcal D}_x^{\beta,\,q})^2$ represents, as equation (3.7)
implies, an unbounded operator on the Hilbert space $L^2(S^1)$ with the scalar
product
$$
\langle g_1,g_2 \rangle =\frac1{2\pi} \int^1_{-1}
\,g_1(x)\,\overline{g_2(x)}\, \widetilde w(x;\beta|\,q)\, dx \,,       \eqno(3.8)
$$
where the weight function $\widetilde w(x;\beta|\,q)$ is defined by (2.3).
In view of (2.2) the polynomials $p_n(x):=d_n^{1/2}(\beta;q)\,C_n(x;\beta|\,q)$,
$n=0,1,2,\cdots$, constitute an orthonormal basis in this space such that
$\left({\mathcal D}_x^{\,q}\right)^2 p_n(x)=  \left(q^{-n/2}+ \beta\,q^{n/2}\right)^2
p_n(x)$. In particular, the ope\-rator $({\mathcal D}_x^{\beta,\,q})^2$ is defined
on the linear span ${\mathcal H}$ of the basis functions $p_n(x)$, which is everywhere
dense in $L^2(S^1)$. We close $\left({\mathcal D}_x^{\beta,\,q}\right)^2$ with respect
to the scalar product (3.8). Since $\left({\mathcal D}_x^{\beta,\,q}\right)^2$ is
diagonal with respect to the orthonormal basis $p_n(x)$, $n=0,1,2,\cdots$, its
closure $\overline{\left({\mathcal D}_x^{\beta,\,q}\right)^2}$ is a self-adjoint
operator, which coincides on ${\mathcal H}$ with $\left({\mathcal D}_x^{\beta,\,q}\right)^2$.
According to the theory of self-adjoint operators (see \cite{AG}, Chapter 6),
we can take a square root of the operator $\overline{\left({\mathcal D}_x^{\beta,\,
q}\right)^2}$. This square root is a self-adjoint operator too and has the same
eigenfunctions as the operator $\overline{\left({\mathcal D}_x^{\beta,\,q}\right)^2}$
does. We denote this operator by $\overline{{\mathcal D}_x^{\beta,\,q}}$. It is evident
that on the subspace ${\mathcal H}$ the operator $\overline{{\mathcal D}_x^{\beta,\,q}}$
coincides with the ${\mathcal D}_x^{\beta,\,q}$. That is, the ${\mathcal D}_x^{\beta,\,q}$
is a well-defined operator on the Hilbert space $L^2(S_1)$ with everywhere dense
subspace of definition. Moreover, according to the definition of a function of a
self-adjoint operator (see \cite{AG}, Chapter 6), we have $\overline{{\mathcal D}_x^{\beta,\,q}}
\,p_n(x)= \left(q^{-n/2}+ \beta\,q^{n/2}\right)\,p_n(x)$. This means that the
continuous $q$-ultraspherical polynomials $C_n(x;\beta|\,q)$ are in fact governed
by a simpler $q$-difference equation,
$$
{\mathcal D}_x^{\beta,\,q}\, C_n(x;\beta|\,q) =  \left( q^{-n/2} +
\beta\,q^{n/2}\right)\,C_n(x;\beta|\,q)\,,                              \eqno(3.9)
$$
which can be regarded as a "square root" of (3.7).

Observe that the $q$-difference operator ${\mathcal D}_x^{\beta,\,q}$ in (3.9) may
be expressed in terms of the Askey-Wilson divided-difference operator $D_q$,
defined in (1.6), as
$$
{\mathcal D}_x^{\beta,\,q}\,=(1+\beta)\,{\mathcal A}_q \,+\,\frac{1-q}{2{\sqrt q}}\,
(1-\beta)\,x\,D_q \,,                                                  \eqno(3.10)
$$
where the ${\mathcal A}_q$ is so-called {\it averaging difference operator}, that is
(see, for example \cite{Ism}),
$$
\left({\mathcal A}_q \,f\right)(x)= \,\frac12\,\left( e^{\,{\rm i}\ln q^{1/2}\,
\partial_{\theta}}+ e^{-{\rm i}\ln q^{1/2}\,\partial_{\theta}}\right)f(x)\,
\equiv\,\cos\left(\ln q^{1/2}\,\partial_{\theta}\right)\,f(x)\,.     \eqno(3.11)
$$

We emphasize that $q$-difference equation (3.8) is consistent with the generating
function
$$
\sum_{n=0}^{\infty}\,t^n \,C_n(x;\beta|\,q)=\,\frac{\left(\beta t\,e^{{\rm i}
\theta},\beta t\,e^{-{\rm i}\theta};q\right)_{\infty}}{\left(t\,e^{{\rm i}
\theta},t\,e^{-{\rm i}\theta};q\right)_{\infty}}                     \eqno(3.12)
$$
for the continuous $q$-ultraspherical polynomials $C_n(x;\beta|\,q)$ (see \cite{GR},
p.169). Indeed, apply the $q$-difference operator ${\mathcal D}_x^{\beta,\,q}$ to both
sides of (3.12) to verify that
$$
\sum_{n=0}^{\infty}\,t^n \,{\mathcal D}_x^{\beta,\,q}\,C_n(x;\beta|\,q)=
{\mathcal D}_x^{\beta,\,q}\,\frac{\left(\beta t\,e^{{\rm i}\theta},
\beta t\,e^{-{\rm i}\theta};q\right)_{\infty}}{\left(t\,e^{{\rm i}\theta},t\,
e^{-{\rm i}\theta};q\right)_{\infty}} $$
$$
= \frac{\left(\,q^{-1/2}\,\beta t\,e^{{\rm i}\theta},\,q^{-1/2}\,
\beta t\,e^{-{\rm i}\theta};q\right)_{\infty}}{\left(\,q^{-1/2}\,t\,
e^{{\rm i}\theta},\,q^{-1/2}\,t\,e^{-{\rm i}\theta};q\right)_{\infty}}
+ \beta \,\frac{\left(\,q^{1/2}\,\beta t\,e^{{\rm i}\theta},\,q^{1/2}\,
\beta t\,e^{-{\rm i}\theta};q\right)_{\infty}}{\left(\,q^{1/2}\,t\,
e^{{\rm i}\theta}, \,q^{1/2}\,t\,e^{-{\rm i}\theta};q\right)_{\infty}}$$
$$
= \sum_{n=0}^{\infty}\,\left(\,q^{-n/2}+ \beta\,q^{n/2}\,\right)\,t^n\,
C_n(x;\beta|\,q)\,.                                              \eqno(3.13)
$$
Equating coefficients of like powers of $t$ on the extremal sides of (3.13),
one completes the another proof of equation (3.9).

As recalled in section 2, if $\beta = q^{\gamma}$, then the $q$-ultraspherical
polynomials $C_n(x;q^{\gamma}|\,q)$ reduce to the Gegenbauer polynomials
$C_n^{(\gamma)}(x)$ in the limit as $q\to 1$. This fact can be also expressed
as the following limit property of the $q$-difference operator ${\mathcal D}_x^{\beta,
\,q}$ in (3.6):
$$
\lim_{q\to 1} \left\{\frac{1}{(\ln q)^2}\,\left[(1 + q^{\gamma})I -
{\mathcal D}_x^{q^{\gamma},\,q}\,\right]\right\} = \frac{1}{4}\,\left[
(1-x^2)\frac{d^2}{dx^2} - (2\gamma + 1)x\frac{d}{dx}\,\right]\,,
$$
where $I$ is the identity operator.

Observe also that the $q$-ultraspherical polynomials $C_n(x;\beta|\,q)$ are
known to possess the simple transformation property
$$
C_n(x;\beta|\,q^{-1})= (\beta q)^n\,C_n(x;\beta^{-1}|\,q)   \eqno(3.14)
$$
with respect to the changes $q\to q^{-1}$ and $\beta \to \beta^{-1}$ (see
\cite{KS}, p.88). It is not hard to check that $q$-difference equation (3.9)
agrees with this property (3.14) since by definition (3.6)
$$
{\mathcal D}_x^{\beta,\,\,q} \equiv \beta\,{\mathcal D}_x^{\beta^{-1},\,\,q^{-1}}\,.
$$

We close this section with the following remark about equation (3.9). Koornwinder
have recently examined raising and lowering relations for the Askey--Wilson
polynomials $p_n(x;a,b,c,d|\,q)$ \cite{K}, which are known to reduce to the
continuous $q$-ultraspherical polynomials $C_n(x;\beta|\,q)$, when one specializes
the parameters $a,b,c,d$ as $a = -\,c= \sqrt{\beta}$ and  $b= -\,d = \sqrt{q\,\beta}$.
So equation (3.9) coincides with "the second order $q$-difference formula" (6.10)
in Koornwinder's paper \cite{K}, upon taking into account that variables $z$ and
$t$ in (6.10) are equal to $e^{{\rm i}\theta}$ and $\beta$, respectively, in our
notations.
\medskip

\setcounter{equation}{0}
\section{Special and limit cases of parameter $\beta$}
\medskip

The $q$-difference equation (3.9) for the $q$-ultraspherical polynomials,
derived in the previous section, does actually contain some special and
limit cases of the parameter $\beta$, which correspond to other well-known
families of $q$-polynomials. We recall (see, for example, \cite{KS}, p.88)
that in the case when $\beta= q^{\alpha + 1/2}$ the $q$-ultraspherical
polynomials $C_n(x;q^{\alpha+1/2}|\,q)$ reduce to (up to a normalization
factor) the continuous $q$-Jacobi polynomials $P_n^{(\alpha,\,\alpha)}(x|\,q)$;
when $\beta= q^{1/2}$ the $C_n(x;q^{1/2}|\,q)$ are related to the continuous
$q$-Legendre polynomials $P_n(x|\,q)$; and when $\beta=q$ the $q$-ultraspherical
polynomials $C_n(x;q|\,q)$ embrace the Chebyshev polynomials of the second
kind $U_n(x)$.

There is also the limit case $\beta \to 1$, which leads to the Chebyshev
polynomials of the first kind $T_n(x)$ in the following way:
$$
\lim_{\beta \to 1}\,\frac{1-q^n}{2(1-\beta)}\,\,C_n(x;\beta|\,q) =
T_n(x)\equiv \cos n\theta\,, \quad n=1,2,3,... \, .
$$
But the point is that $q$-difference equation (3.9) in this limit reduces to
the difference equation
$$
\left[\,e^{{\rm i}\ln q^{1/2}\partial_{\theta}} + e^{-{\rm i}\ln q^{1/2}
\partial_{\theta}}\right]\,T_n(x) = \left(q^{n/2} + q^{-n/2}\right)\,T_n(x)\,,\eqno(4.1)
$$
although we all know well that the Chebyshev polynomials of the first kind $T_n(x)$
satisfy the second-order differential equation
$$
\left[\,(1-x^2)\,\frac{d^2}{dx^2} - \,x\,\frac{d}{dx}+ n^2\,\right]T_n(x)=0\,.
$$
Nevertheless, there is no contradiction here since one readily verifies that the
Chebyshev polynomials of the first kind $T_n(x)=\cos{n\theta}$,\, $n=0,1,2,...$,
do satisfy difference equation (4.1) as well.


\setcounter{equation}{0}
\section{Concluding remarks}
\medskip

To summarize, we have proved that the conventional $q$-difference equation
(2.4) of Sturm-Liouville type for the continuous $q$-ultraspherical
polynomials $C_n(x;\beta|\,q)$ of Rogers admits factorization of the form
(3.9). The special case of the $C_n(x;\beta|\,q)$ with the vanishing
parameter $\beta$ is known to correspond to the continuous $q$-Hermite
polynomials $H_n(x|\,q)$. The above-presented formulas in this case when
$\beta=0$ are in accord with that obtained by M.Atakishiyev and A.Klimyk
in \cite{AK}. So it would be of considerable interest to explore now
whether the situation here described obtains for other families of
orthogonal polynomials on higher levels in the Askey $q$-scheme \cite{KS}.
Work on clarifying this point is in progress.

We are grateful to N. Atakishiyev and E. Godoy for encouraging our interest in this
problem and helpful discussions.

\medskip
\section*{Acknowledgements}
\medskip

MKA would like to thank the Departamento de Matem\'atica Aplicada II,
Universidade de Vigo, Spain for their hospitality during her visit in
April, 2007 when the main part of this research was carried out.

\end{document}